\documentclass[12pt]{article}
\usepackage{amssymb}
\usepackage{theorem}

\newtheorem{theorem}{Theorem}[section]
\newtheorem{proposition}[theorem]{Proposition}
\newtheorem{corollary}[theorem]{Corollary}
\newtheorem{lemma}[theorem]{Lemma}
{\theorembodyfont{\rmfamily}
\newtheorem{definition}[theorem]{Definition}

\newtheorem{remark}[theorem]{Remark}

}

\include{epsf}

\begin{document}

\def\CP{\mathbb{C}{\rm P}}
\def\tr{{\rm tr}\,}
\def\endproof{{$\Box$}}

\title{Complex surfaces with $\rm CAT(0)$ metrics}

\author{Dmitri Panov
\footnote{Supported by a Royal Society University Research Fellowship}}

\maketitle

\begin{abstract}

We study complex surfaces with
locally $\rm CAT(0)$ polyhedral
K\"ahler metrics and construct such metrics
on $\mathbb CP^2$ with various orbifold structures.
In particular, in relation to questions of Gromov and Davis-Moussong
we construct such  metrics
on a compact quotient of the two-dimensional unit complex ball.
In the course of the proof of these results we give
criteria for Sasakian $3$-manifolds
to be globally $\rm CAT(1)$.
We show further that for certain Kummer coverings of $\mathbb CP^2$
of sufficiently high degree
their desingularizations are of type $K(\pi,1)$.

\end{abstract}

\tableofcontents

\section{Introduction}

In this article we  study curvature properties of
{\it polyhedral K\"ahler manifolds} introduced in \cite{P}.
Before recalling the definition of polyhedral K\"ahler
manifolds it is worth to describe
two types of examples that are rather familiar.
The first class are quotients of  complex tori with a flat K\"ahler
metric by a finite group of isometries. The second class are
ramified covers of $\mathbb C^n$ branched along
a collection of complex hyperplanes (these examples are not compact).
In both cases we have a complex manifold (smooth or singular)
with a flat metric with singularities. Metric singularities
of complex co-dimension $1$ occur along flat geodesic strata.
In the first
case the {\it conical angles} in the normal direction to such
strata are $\frac{2\pi}{k}$ while in the second they are $2\pi k$,
$k\in \mathbb N$.
Polyhedral K\"ahler manifolds have similar structure, but
the conical angle at geodesic divisors can be arbitrary.

Now, recall the definition.
A {\it polyhedral K\"ahler manifold} is a $PL$
manifold $M^{2n}$
with a flat  metric with singularities and
a compatible parallel complex structure defined outside of metric
singularities. The metric should satisfy two conditions.
1) $M^{2n}$
should admit a simplicial decomposition
such that the metric restricts to Euclidean metric on each
simplex. As a consequence the singularities of the metric
represent a  union of several $2n-2$--dimensional simplexes.
2) Every  $2n-2$-dimensional simplex contained in the singular locus
of the metric is holomorphic
with respect to the complex structure on each top dimensional simplex
containing it. We use abbreviation $PK$ to denote polyhedral
K\"ahler manifolds and metrics on them.

It was proven in \cite{P} that $4$-dimensional $PK$ manifolds
are smooth complex surfaces and the metric singularities
form a union of complex curves on these surfaces.
For the purposes of this article we widen the class of $PK$ manifolds in the
$4$-dimensional case and consider $PK$
metrics defined on complex surfaces (of complex dimension two) with isolated complex
singularities.

We give criteria on the singularities of $PK$ metrics
implying that the metric is locally $\rm CAT(0)$ (see Definition \ref{def:negative}),
and use them to construct
several families of locally $\rm CAT(0)$ orbifold metrics on $\mathbb CP^2$.
In particular we produce a locally
$\rm CAT(0)$ polyhedral metric on a compact complex hyperbolic $2$-ball
quotient (see Section \ref{threearrangements}). This is related to questions raised by Gromov (\cite{Gr2} , Section 7.A.IV, p. 180),
and Davis with Moussong (see discussion in \cite{DM} Question 2, Section 7.3).
Further we construct several
families of complex surfaces of type $K(\pi,1)$.

\subsection{Criteria for being locally $\rm CAT(0)$}
A metric on a polyhedral K\"ahler manifold is flat
outside of singularities, so to check that it is
locally $\rm CAT(0)$ we need to study metric singularities.
On a $4$-dimensional  $PK$ manifold metric singularities
occur in real codimension  $2$ (along complex curves)
and in codimension $4$ (where these
curves meet). The
first case is easy, the metric is locally a direct product
of a $2$-dimensional cone with $\mathbb R^2$; it is locally
$\rm CAT(0)$ if an only if the $2$-cone has
conical angle larger than $2\pi$. The real task central for this
article is to understand
the nature of locally $\rm CAT(0)$ condition at metric
singularities of codimension $4$.

For every point $x$ of a $PK$
$4$-manifold we denote by $C_K^4(x)$ its
{\it tangent cone}, i.e., the cone that has a neighborhood of its origin
isometric to a neighborhood of $x$. By $S^3_K(x)$ we denote
the unit sphere in the tangent cone.
A metric on a $PK$ surface
is locally $\rm CAT(0)$ if and only if for every point $x$ the sphere
$S^3_K(x)$ is $\rm CAT(1)$. By Corollary \ref{prop:cones}
this happen if and only if
any closed geodesics in $S^3_K(x)$ is not shorter than
$2\pi$, and the conical angle at each complex
curve  from the singular locus is
greater than $2\pi$.

Recall \cite{P} that for any point $x$ on a $PK$ $4$-manifold the
cone $C_K^4(x)$ admits
an isometric holomorphic action by $\mathbb R^1$, given by
the field $ir\frac{\partial}{\partial r}$ (here
$r\frac{\partial}{\partial r}$ is the usual dilatation field that
exists on every polyhedral cone).
If $C_K^4(x)$ is not a direct product of two cones
of angles $2\pi \beta_1$ and  $2\pi \beta_2$ with
$\frac{\beta_1}{\beta_2}\notin \mathbb Q$,
then  the isometric
action of $\mathbb R^1$ factors through an $S^1$-action.
In this article we concentrate
on the cones where this $S^1$ action on $S_K^3(x)$ is free,
and we call such cones {\it regular spherical cones}.
It was shown in \cite{P} (Theorem 1.7)
that for a regular spherical cone there exists
a canonical (up to scaling) biholomorphism $C_K^4(x)\to \mathbb C^2$
that sends the singular locus of $C_K^4(x)$ to a collection
of lines passing through the origin of $\mathbb C^2$.
It will also be important for us to consider ramified
covers of regular spherical cones with branching at
the lines where the metric has singularities. These are
polyhedral cones not biholomorphic to $\mathbb C^2$
and we call them {\it regular} if the natural action of $S^1$
on them is free.

The next two theorems provide two main classes of
locally $\rm CAT(0)$ metric singularities that appear in this article.

\begin{theorem}\label{CAT1S3}
Let $S_K^3$ be the unit sphere of a  regular spherical $PK$
cone $C_K^4$.
If the quotient  $2$-sphere $S_K^3/S^1$ is
$\rm CAT(4)$ then the sphere  $S_K^3$
is $\rm CAT(1)$ and consequently the
cone $C_K^4$ is locally $\rm CAT(0)$.
\end{theorem}

Theorem \ref{CAT1S3}
is related to results obtained \cite{CD}. Namely,
Theorem 9.1 from \cite{CD}
gives a {\it necessary} and {\it sufficient}
criterion for polyhedral metrics on $\mathbb R^4$
whose quotient by a complex reflection
group is isometric to the flat $\mathbb R^4$,
to be $\rm CAT(0)$.
Theorem \ref{CAT1S3} extends
the {\it sufficiency} conditions from \cite{CD}
to the whole class of spherical $PK$ cones.
We prove as well a smooth analogue of Theorem \ref{CAT1S3}
for {\it Sasakian} metrics on $S^3$, Theorem \ref{smoothCAT1S3}.

Next theorem explains that locally $\rm CAT(0)$
cones can be constructed by taking ramified covers
of non-negatively curved cones.

\begin{theorem}\label{COVERCAT1}
Let $S_K^3$ be the unit sphere of a regular spherical  $PK$ cone
$C_K^4$. Let $s_1,...,s_k$, $k\ge 3$ be the singular
circles of the metric on $S_K^3$. Suppose that the conical angles along
the circles $s_1,...,s_k$ are less than $2\pi$. Then for any
sufficiently large integers $n_1,...,n_k>0$ there exists a $3$-manifold $M^3$
of type $\rm CAT(1)$ with an isometric action of a finite group $G$
such that the quotient $M/G$ is isometric to $S_K^3$ and has branching
of order $n_i$ along $s_i$.
\end{theorem}

\subsection{Constructions of complex surfaces and orbifolds}

{\bf Orbifold $\rm CAT(0)$ metrics on $\mathbb CP^2$.}
In order to construct orbifold $\rm CAT(0)$ metrics on
$\mathbb CP^2$ we  use line arrangements
denoted by $A_1(6)$, $A_1(7)$ and $A_3^0(3)$
in \cite{G} and \cite{H} (for the definition see Section \ref{threearrangements}). These arrangements satisfy two important
properties. First, they are  singular loci of $PK$ metrics on $\mathbb CP^2$
with conical angles less than $2\pi$ (\cite{P}).
Second, the multiple points of these arrangements have multiplicity
at most three. Each of these arrangements gives rise to an infinite series of
locally $\rm CAT(0)$ orbifold structures on $\mathbb CP^2$,
where to each line we associate a certain orbifold multiplicity
$d$ so that the orbifold stabilizer at non-multiple
point of the line is $\mathbb Z_d$. In Section  \ref{orbiresults}
we explain further how the orbifold structure is defined at multiple points.

\begin{theorem}\label{CPorbi}
Let $(L_1,...,L_n)$ be one of the arrangements
$A_1(6)$, $A_1(7)$ or $A^0_3(3)$. Associate integers
$d_i>1$ to lines $L_i$
in such a way that at every triple point of the
arrangement incoming numbers $d_i, d_j, d_k$
satisfy $\frac{1}{d_i}+\frac{1}{d_j}+\frac{1}{d_k}>1$.
Then for each orbifold structure $(L_1,d_1;...;L_n,d_n)$ on
$\mathbb CP^2$ there exits a $PK$ metric that is locally
$\rm CAT(0)$ with respect to the orbifold structure.
\end{theorem}

Orbifolds in Theorem \ref{CPorbi} have contractible universal
covers. In several cases the universal cover  is bi-holomorphic to $\mathbb C^2$,
polidisk, or unit complex ball. Some of these cases have appeared in \cite{U}, and all such
orbifods are representable as finite quotients
of smooth complex surfaces. However for majority
of orbifolds in Theorem \ref{CPorbi} the universal cover is not of one of these three types.
One can expect that such orbifolds are finite quotients
of negatively curved surfaces, analogous to
surfaces of Mostow-Siu  \cite{MS}. See also Remark \ref{cat(-1)}
where we explain how one can find an orbifold $\rm CAT(-1)$ metric
in the case of arrangement $A^0_3(3)$.

{\bf Aspherical complex surfaces.}
Using Kummer covers of $\mathbb CP^2$ and
Theorem \ref{COVERCAT1} we produce several series
of aspherical complex surfaces of type $K(\pi,1)$.
Recall that Kummer covers are
associated to line arrangements $(L_{1},...,L_{k})$
on $\mathbb CP^2$ and
an integer $n>0$. One  considers the extension of the field
of rational functions on $\mathbb{C}P^{2}$ by functions
$(\frac{L_{i}}{L_{j}})^{\frac{1}{n}}$. This extension
defines a Galois cover
$S(n,L_{1},...,L_{k})\to \mathbb{C}P^{2}$ of degree $n^{k-1}$.

\begin{theorem}\label{ramCP2}
Consider a $PK$ metric on $\mathbb CP^2$ that has conical angles
less than $2\pi$ at the line arrangement $L_{1},...,L_{k}$.
There exists $N$ such that for every
$n\ge N$ the blow up of $S(n,L_{1},...,L_{k})$
at its complex singularities is an aspherical complex
surface.
\end{theorem}

The fundamental groups of these surfaces are quite non-trivial
infinite quotients of the fundamental group of the complement
to arrangement $(L_i)$ in $\mathbb CP^2$.
The proof of the theorem relies on the following general statement.

\begin{theorem}\label{kpi1}
Let $S$ be complex surface with isolated complex singularities. If
$S$ admits a locally $\rm CAT(0)$ metric then any blow up of
$S$  that does not contain rational exceptional
curves  is of type $K(\pi,1)$.
\end{theorem}

We note finally that
Kummer covers were extensively studied starting
form the work of Hirzebruch \cite{H},
who used them to construct
complex ball quotients that are blow ups of $\mathbb CP^2$
at several points.

{\bf Acknowledgments.} This article is a substantial revision of
the final part of my PhD written at \'Ecole
Polythechnique under supervision of Maxim Kontsevich
and defended in 2005. I am very indebted to Maxim
for his help and advise. Numerous discussions with
Anton Petrunin were essential for this work.
I would like to thank  Robert Bryant, Martin Deraux,
Misha Gromov, Christophe Margerin, and Dimitri Zvonkine for
discussions and help. I would also like to thank
Tadeusz Januszkiewicz for the encouragement to finish this paper, Robert Bell for the reference \cite{Bo}, and an anonymous referee for helpful suggestions.
This work was supported by EPSRC grant EP/E044859/1.

\section{\label{Catrecall}
Generalities about $\rm CAT(\kappa)$ spaces
and proof of Theorem \ref{kpi1} }

In this section we mainly recall definitions and theorems about
$\rm CAT (\kappa)$ spaces.
The material is borrowed from  Chapter 10 of
\cite{GH} and Chapter 4 of \cite{Gr1} unless
a different source is specified.

Denote by $M_{\kappa}^k$  a
complete simply connected manifold of constant
curvature $\kappa$, i.e., the hyperbolic $k$-space for $\kappa<0$,
the flat
$\mathbb {R}^k$ for $\kappa=0$ and the $k$-sphere (of radius
$\frac{1}{\sqrt{\kappa}}$) for $\kappa>0$.

\begin{definition}\label{catdefinition}
Let $X$ be a complete geodesic space.
A {\it triangle} $\Delta$
in $X$ is a union of three geodesic segments that join three points of $X$.
The {\it comparison triangle} of $\Delta$ in $M_{\kappa}^2$ is the triangle
$\tilde{\Delta}$ with the edges
of the same length as $\Delta$. We say that
$\Delta$ satisfies the $\rm CAT(\kappa)$ inequality if
for all $x,y\subset \Delta$ and the corresponding points
$\tilde{x}, \tilde{y}\subset\tilde{\Delta}$  the following inequality holds:
$$d(x,y)\le d(\tilde{x},\tilde{y}).$$
A space $X$ is called a $\rm CAT(\kappa)$
space if every triangle in
$X$ of perimeter less than $\frac{2\pi}{\sqrt{\kappa}}$ satisfies
the $\rm CAT(\kappa)$ inequality (if $\kappa\le 0$
the inequality should be satisfied for all triangles).
\end{definition}

\begin{definition}\label{def:negative}
A space $X$ has curvature $K_X\le\kappa $ if
for every point $p\in X$ there is a neighborhood $U_p$ such that
every triangle in $U_p$ of perimeter less than
$\frac{2\pi}{\sqrt{\kappa}}$ satisfies the $\rm CAT(\kappa)$
inequality (if $\kappa\le 0$ we  consider all triangles).
If $K_X\le \kappa$ we also say that $X$ is
{\it locally $\rm CAT(\kappa)$}.
\end{definition}

\begin{theorem}\label{thm:CAT}
Suppose that $K_X\le \kappa$
 and every closed geodesic in $X$ has length at least
$\frac{2\pi}{\sqrt{\kappa}}$.
Then $X$ is a $\rm CAT (\kappa)$ space.
\end{theorem}

\begin{theorem}\label{thm:KP1}
A simply connected geodesic space $X$ with $K_{X}\le 0$
is contractible, every two points on it are joined by a unique geodesic.
\end{theorem}

{\bf Proof of Theorem \ref{kpi1}.}
Let $\widetilde S$ be the universal cover of
$S$, and $\widetilde S_{\circ}$
the  blow up of $\widetilde S$ that covers
the chosen blow up of $S$.  To prove the theorem
it is sufficient to show that $\widetilde S_{\circ}$
is $K(\pi,1)$.
Since
$\widetilde S_{\circ}$ is a blow up of $\widetilde S$
in a discrete subset of complex singularities,
any homotopy of a sphere in
$\widetilde S_{\circ}$ can touch only finite number of exceptional curves.
Hence it is enough to
prove that the blow up of $\widetilde S$ at any finite
subset of complex singularities, say $x_1,...,x_n$ is $K(\pi,1)$.

Fix a point $x$ of $\widetilde S$.
By Theorem \ref{thm:KP1} for any $1\le i\le n$
there exists a unique
geodesic $\gamma_i$ on $\widetilde S$ that joins  $x$ and $x_i$.
The union of the geodesics $\gamma_i$ forms a tree  on
$\widetilde S$, we denote this tree by $T$
and denote by $T_{\varepsilon}$ the $\varepsilon$-neighborhood of $T$.
Since $\widetilde S$ is $\rm CAT(0)$  every geodesic ray starting
at $x$ intersects the boundary $\partial T_{\varepsilon}$
exactly at one point.  Indeed, in any $\rm CAT(0)$ space for any
geodesic ray $r(t)$ and a geodesic segment $s(t)$ sharing
common vertex $x=r(0)=s(0)$ the distance from $r(t)$
to the segment $s$ is a function increasing with $t$. It follows
that $T_{\varepsilon}$ is contractible,
and as well $\widetilde S$ can be contracted on
$T_{\varepsilon}$ by
a homotopy identical on $T_{\varepsilon}$ and preserving geodesics
rays starting at $x$.

Notice finally that by the assumption of the theorem
the blow up of $T_{\varepsilon}$ at
$x_1,...,x_n$ is homotopic to a bouquet of $n$
spaces of type $K(\pi,1)$. Indeed,  the blow
up of $T_{\varepsilon}$ is homotopic
to the bouquet of $\varepsilon$-neighborhoods
$U_{\varepsilon}( x_i)$ of $x_i$
blown up at points $x_i$.
At the same time  every $U_{\varepsilon}(x_i)$
blown up at $x_i$ can be contracted on
its exceptional curve, which is a union of several
curves of genus higher than $0$ with several points identified.
Such spaces  are $K(\pi,1)$.

\hfill $\square$

{\bf Polyhedra of curvature $\kappa\le K$.}
We use the terminology from \cite{Gr1}.
A {\it straight simplex} in $M^k_{\kappa}$
is a compact intersection
of $k+1$ half-spaces in $M^k_{\kappa}$ in generic position.
Let $X$ be a finite-dimensional simplicial complex
with a metric such that each simplex is isometric to
a straight simplex in $M_{\kappa}^k$ for some $k$.
Such a space $X$ is called $(M,{\kappa})$ space. To each
$k$-simplex $\Delta$ on $X$ we can associate its {\it link}
$L_{\Delta}$. This is an $(M,1)$ space that can be defined
as the set of unit vectors orthogonal to $\Delta$ at an
interior point. Here to each simplex $\Delta'$ containing
$\Delta$ corresponds a simplex in $L(\Delta)$
of curvature $1$ and of dimension $dim(\Delta')-k-1$.
The following theorem is stated, for example
on page 2 in \cite{Bo}.

\begin{theorem}\label{thm:negativepolyhedra}
Let $X$ be a $(M,\kappa_0)$ space.
Then $K_X\le \kappa$
if and only $\kappa_0\le \kappa$ and
for any $\Delta$ in $X$ the link $L_{\Delta}$
does not contain a closed geodesic of length
less than $2\pi$.
\end{theorem}

Of course, a $PK$ manifold is an $(M,0)$ space.
In order to apply Theorem \ref{thm:negativepolyhedra}
to a $4$-dimensional $PK$ cone
we need to verify the $\rm CAT(1)$
condition on the link at its origin and as well
at all complex curves forming the singular locus of the cone.
For such curves the link is just the circle of length
equal to the conical angle at the curve. This is summarized
in the following corollary:

\begin{corollary} \label{prop:cones}
A  $PK$ cone $C_K^4$ is locally $\rm CAT(0)$ if and only if the conical angle
at any curve in its singular locus is greater than $2\pi$ and
closed geodesics on its unitary sphere $S_K^3$ have  length
at least $2\pi$.
\end{corollary}

The following result on $\rm CAT(\kappa)$
spaces can be found in  \cite{AKP}.

\begin{theorem} \label{convergence}
Let $X_n$ be a sequence of $\rm CAT(\kappa)$ spaces that
converges in Gromov-Hausdorff topology to $X$.
Then $X$ is also a $\rm CAT(\kappa)$ space.
\end{theorem}

Finally we quote a comparison theorem from \cite{AB}
on curves with  bounded geodesic curvature on $\rm CAT(\kappa)$
spaces. The definition of geodesic curvature for curves in
$\rm CAT(\kappa)$ spaces is introduced on page 70 in \cite{AB},
it is called
{\it arch-chord curvature}. We don't reproduce the
definition here since it is a bit technical, while
the comparison theorem will be applied at worst to
piecewise smooth curves on real surfaces with conical points.
By a {\it loop} in a $\rm CAT(\kappa)$ space we mean
a rectifiable curve $\gamma$ that admits a
parameterization $\phi[0,1]\to \gamma$ with
$\phi(0)=\phi(1)$.   The following statement is a special case of
Theorem 1.1 \cite{AB} that treats as well curves with
$\phi(0)\ne \phi(1)$.

\begin{theorem}\label{piloop}
 Let $X_{\kappa}$ be $\rm CAT(\kappa)$ space
and let $\gamma$ a loop in it of geodesic curvature
at most $k$ (the curvature need not be defined at the base
point of $\gamma$). Then $\gamma$ is no shorter than
the complete circle of curvature $k$ on $M^2_{\kappa}$.
\end{theorem}

\section{$\rm CAT(1)$ property of unit spheres of $PK$ cones}

The goal of this section is to prove Theorems \ref{CAT1S3}
and \ref{COVERCAT1}. Let us outline briefly the strategy
of the proofs.
In both theorems we need to show that
on the unit sphere $S_K^3$ of a certain regular
$PK$ cone all closed geodesics
have length at least $2\pi$.
The first step is to  project geodesics in $S_K^3$ to
the  quotient $S_K^3/S^1$ and notice that the projection is
a curve of constant curvature on the complement to
conical points of $S_K^3/S^1$ (Lemma \ref{geonSK3}).
One can show then that in the case when
$S_K^3$ is topologically $S^3$ and
conical points on $S_K^3/S^1$
have angles greater than $2\pi$, every closed geodesic on
$S_K^3$ projects to a curve with a self-intersection
on $S_K^3/S^1$ (Proposition \ref{selfint}).
This observation joined with Theorem \ref{piloop}
settles Theorem \ref{CAT1S3}. Theorem  \ref{COVERCAT1}
is deduced from  its
$2$-dimensional analogue (Corollary \ref{good1} ) joined
with Theorem \ref{piloop}.

In addition to Theorem \ref{CAT1S3} we prove its smooth
version (Theorem \ref{smoothCAT1S3}) for Sasakian $3$-manifolds.
Both proofs follow the same path and for the sake of
exposition we give them simultaneously,
giving all the details for the smooth case
and explaining how to alter them for the original case of
Theorem \ref{CAT1S3}. In principle, using carefully Theorem
\ref{convergence}
one could also show that the statements of
Theorem \ref{CAT1S3} and its smooth version \ref{smoothCAT1S3}
are equivalent. Indeed every singular metric on $S^3$ from \ref{CAT1S3}
can be Gromov-Hausdorff approximated by Sasakian
metrics from \ref{smoothCAT1S3} and vice versa.
We chose not to do this to avoid technical issues
that arise in this approach.

\subsection{Geodesics of unit spheres of $PK$ singularities}

\begin{lemma}\label{geonSK3}
Let $\gamma$ be a geodesic on the unit sphere $S_K^3$
of a regular $PK$ cone.
Then the angle $\angle\gamma'$  between $\gamma$
and the $S^1$-fibration is constant.
The lengths of $\gamma$ and it projection $p(\gamma)$ to $S^3_K/S^1$
are related by the formula
$l(\gamma)\mathrm{sin}(\angle\gamma')=l(p(\gamma))$.
Moreover $p(\gamma)$ is a piecewise smooth curve of constant
curvature
$2\frac{\mathrm{cos}(\angle \gamma')}{\mathrm{sin}(\angle \gamma')}$ that
can have singularities only at conical points of $S_K^3/S^1$.

2) Let $x\in p(\gamma)$ be a conical point on $S^3_K/S^1$.
Then locally $p(\gamma)$ divides a neighborhood of $x$ into two
sectors with angles at $x$ greater or equal to $\pi$.
\end{lemma}

{\bf Proof.} All these statements are obvious in the
case when $S_K^3$ is the unit sphere
of the flat $\mathbb C^2$, with the standard $S^1$ action
defining the Hopf fibration on $S_K^3$.
In this case the quotient $S^3_K/S^1$
is the sphere of curvature $4$.

1) In the general case any closed geodesic $\gamma$
on the unit sphere $S_K^3$ of a regular $PK$ cone
is composed of smooth geodesic segments
that join singular fibers of the $S^1$-fibration. The function
$\angle \gamma'$ is constant on each smooth segment of $\gamma$ because
the $S^1$-fibration on the complement
to singularities in $S_K^3$ is
locally isometric
to the Hopf fibration on the standard sphere $S^3$.

Let us show that the value of $\angle \gamma'$
is the same on all smooth segments of $\gamma$.
Denote by $T^2=S^1\times \gamma$   the torus immersed in
$S_K^3$ that is the union of $S^1$ fibers intersecting
$\gamma$. The restriction metric on $T^2$ is flat,
$\gamma$ is embedded in  $T^2$ as a geodesic and so it intersects
the flat foliation given by  $S^1$-fibers under constant angle.
Hence we have $\angle \gamma'=\mathrm{const}$ and the further
claims of the lemma hold since they hold for the standard unit sphere.

2) If one of two sectors of $p(\gamma)$ at $x$
had angle smaller than $\pi$  one would be able to smoothen $\gamma$
in a neighborhood of $p^{-1}x$ to make it shorter,
which is impossible since $\gamma$
is locally minimizing.

\hfill $\square$

\subsection{Geodesics on Sasakian $3$-manifolds}
In this section we recall the definition of
regular {\it Sasakian}  $3$-manifolds and study geodesics on them.
Let  $M^3\to M^2$ be an $S^1$-bundle over a real surface $M^2$.

\begin{definition}A metric $g$  on $M^3$ is called regular
{\it Sasakian} if it satisfies
the following  conditions:

1) There is an isometric action of $S^1$ on $M^3$ that preserves the fibers

2) Let $\omega_0$ be the $1$-form of norm $1$ with zero distribution
orthogonal to the $S^1$-fibers. Then $d\omega_0$
equals $-2d{\rm S}$, where $d{\rm S}$ is the area form on $M^3/S^1$.

\end{definition}

Regular Sasakian manifolds are special case of Sasakian \cite{BG},
that don't necessarily  admit a structure of $S^1$-fibration.
The simplest example of a Sasakian $3$-manifold is the unit sphere
where the $S^1$-fibration is the Hopf fibration. Unit spheres
or $PK$ cones can be considered as examples of Sasakian manifolds with
singular metric.

The following is the main statement of this section, I would like
to thank Robert Bryant for the help with its proof.

\begin{proposition}\label{Bryant}
Let $\gamma$ be a geodesic
on a regular Sasakian manifold $M^3$.
The angle between $\gamma$ and the $S^1$-fibration is constant and
the projection of $\gamma$ to $M^2$ is a smooth curve
of constant curvature equal to $2\omega_0$.
\end{proposition}

{\bf Proof.} Let $\omega_0,\omega_1,\omega_2$
be  $S^1$-invariant  $1$-forms on
$M^3$ forming  an orthogonal base of $T^*M^3$, and
such that the restriction of $\omega_1$ and $\omega_2$
on $S^1$-fibers are zero. Let $e_0,e_1,e_2$ be the dual
vector fields on $M^3$. Then $[e_0,e_1]=[e_0,e_2]=0$, and
$[e_1,e_2]=2e_0$.  Plugging this in
the formula for the exterior derivative of $1$-forms
$d\omega(X,Y)=L_X\omega(Y)-L_Y\omega(X)-\omega([X,Y])$
we get:

\begin{equation}\label{differentials}
d\omega_0=-2\omega_1\wedge\omega_2, \;\;\;
d\omega_1=-\theta_{12}\wedge\omega_2,\;\;\;
d\omega_2=\theta_{12}\wedge\omega_1,
\end{equation}
where $\theta_{12}$ is a uniquely defined $1$-form,
that is  a (point-wise) linear combination of $w_1$ and $w_2$.

Let $p_0,p_1,p_2$ be the coordinates on fibers of $T^*M^3$ dual to
$\omega_0,\omega_1,\omega_2$.
Then the canonical $1$-form $\alpha$ on $T^*M^3$  is
the following:

\begin{equation}\label{canform}
\alpha=p_0\omega_0+p_1\omega_1+p_2\omega_2.
\end{equation}

Using  formula (\ref{differentials}) and (\ref{canform})
we get the following expression for the
symplectic form on $T^*M^3$:

\begin{equation}\label{dalpha}
d\alpha=dp_0\wedge\omega_0+
(dp_1+p_2\theta_{12}-p_0\omega_2)\wedge\omega_1+
(dp_2-p_1\theta_{12}+p_0\omega_1)\wedge\omega_2.
\end{equation}

Denote by $X$ the geodesic flow on $T^*M^3$. From Lee formula it follows

\begin{equation}\label{lider}
i(X)d\alpha=-\alpha=-p_0dp_0-p_1dp_1-p_2dp_2.
\end{equation}

Using $\omega_iX=p_i$ we deduce form (\ref{dalpha}) and (\ref{lider})
the  following system of equations
$$dp_0 X=0,\;\;\;(dp_1+p_2\theta_{12}-p_0\omega_2)X=0,
\;\;\;(dp_2-p_1\theta_{12}+p_0\omega_1)X=0.
$$

Here the first equation states that any geodesic on $M^3$
has fixed angle with the $S^1$ fibration, namely the
quantity $p_0(X)$ is constant along the trajectories
of the geodesic flow. On the over hand, if we fix $p_0(X)=c$
the next two equations can be interpreted
as the equations for the twisted
geodesic (magnetic) flow on $T^*(M^3/S^1)$, whose trajectories
project to curves of geodesic curvature $2c$ on $M^3/S^1$.
In particular by choosing $c=0$ we get the geodesic flow.
To go from $M^3$ to $M^3/S^1$ we use the $S^1$-invariance of
all equations.

\hfill $\square$

\begin{lemma} \label{2area}
Consider a regular Sasakian metric on $S^3$ and
let $S^1\hookrightarrow S^3\to S^2$ be the corresponding fibration.
Then the equally holds $l(S^1)=2{\rm Area}S^2$.
The same relation holds for the unit
sphere of a regular spherical $PK$ cone.
\end{lemma}

{\bf Proof.} In the Sasakian case applying Stokes formula to a disk
in $S^3$  whose boundary is a fiber $S^1$
and using the definition of Sasakian metric we get
$\int_{S^1}\omega_0=-\int d\omega_0=2{\rm Area(S^2)}$.
The second statement is a special case of Theorem 1.9 \cite{P}.

\hfill $\square$

Finally we give one standard proposition.

\begin{proposition}\label{threestandard}
The sectional curvature of a regular Sasakian $3$-manifold $M^3$
is bounded from above by $1$ if $K(M^3/S^1)\le 4$.

\end{proposition}

{\bf Proof.} This can be found, for example in Section 7 \cite{BG}.
The sectional curvature of each $2$-plane containing the vertical
direction equals $1$, and the curvature of a horizontal $2$-plane
equals $K(M^3/S^1)-3$ (this can also be deduced
from O'Neill's formula).

\hfill $\square$

\subsection {$S^3$ is $\rm CAT(1)$ if $S^2$ is $\rm CAT(4)$}

In this subsection we prove Theorem \ref{CAT1S3} and its smooth analogue.

\begin{proposition}\label{twoineq}
1) Let $\gamma$ be a closed geodesic in a regular Sasakian manifold
$M^3$, whose projection to $M^3/S^1$
bounds a topological disk $\Omega_{\gamma}$.
Stokes formula on $M^3$ applied to $(\omega_0,\gamma)$
and, Gauss-Bonnet applied to $\gamma$ on $M^3/S^1$
give us:

\begin{equation}\label{modl(S)}
\int_{\gamma}w_0\equiv
-2Area (\Omega_{\gamma})\;\; {\rm mod}\;\; l(S^1),
\end{equation}
\begin{equation}\label{gaussbonnet}
\int_{\gamma}w_0=
\frac{1}{2}(2\pi-\int_{\Omega_{\gamma}}K\;{\rm dS}).
\end{equation}
2) Suitably interpreted Equalities (\ref{modl(S)}) and
(\ref{gaussbonnet} )   hold as well
when $M^3$ is a unit sphere of a regular polyhedral K\"ahler cone.
\end{proposition}
{\bf Proof.} 1) To prove Equality (\ref{modl(S)}) consider
a section $s$ of the $S^1$ bundle over $\Omega_{\gamma}$.
Let $\partial s$ be the restriction of this section over
$\partial \Omega_{\gamma}$. Using the
definition of Sasakian metric and applying Stokes formula
we get
$$-2Area (\Omega_{\gamma})=\int_S dw_0=
\int_{\partial S}w_0\equiv\int_{\gamma}w_0 \;\; mod\;\; l(S^1).$$
Equality (\ref{gaussbonnet} )
follows from Proposition \ref{Bryant} together
with Gauss-Bonnet formula.

2) In the case when $M^3$ is the unit sphere of a
$PK$ cone,
Equation (\ref{modl(S)}) does not require any interpretation
and is proved as above,
one just needs to use Lemma \ref{geonSK3} instead of
Proposition \ref{Bryant}. To interpret correctly
(\ref{gaussbonnet}) we first notice that the curvature
$K$ on $M^3/S^1$ is the sum of the constant $4$ with delta-functions,
supported at conical points.
A conical point $x$ of
conical angle $2\pi\alpha(x)$ supports the delta function
$\delta(x)=2\pi(1-\alpha(x))\delta$.
In order to interpret correctly
the integral $\int_{\Omega_{\gamma}}K\;{\rm dS}$ we should
specify the contribution to it of conical points lying on the
boundary of $\Omega_{\gamma}$. Let $y$ be such a point
and let $\pi\alpha_{in}(y)$ be the angle between two arcs of
$\gamma$ measured inside of $\Omega_{\gamma}$, and
$\pi\alpha_{out}(y)$ be the angle measured in the complement
$(M^3/S^1)\setminus \Omega_{\gamma}$. Of course, the sum
$\pi\alpha_{in}(y)+\pi\alpha_{out}(y)$ equals the total conical
angle $2\pi\alpha(y)$  at $y$, and $\alpha_{in}\ge 1$,
$\alpha_{out}\ge 1$ since $\gamma$ is a geodesic.
We split the delta function $\delta(y)$ as a sum of
$\delta_{in}(y)=\pi(1-\alpha_{in}(y))\delta$ and
$\delta_{out}(y)=\pi(1-\alpha_{out}(y))\delta$.
Assigning to each conical point $y$ on the boundary of $\Omega_{\gamma}$
the delta function $\delta_{in}(y)$ one gets the correct version
of formula  (\ref{gaussbonnet}).

\hfill $\square$

\begin{proposition}\label{selfint}
Consider a regular Sasakian metric on $S^3$ and let  $S^3\to S^2$ be
the corresponding  $S^1$-fibration.
Suppose that $K(S^2)\le 4$. Then for any  closed
geodesic $\gamma$ on $S^3$ its projection to $S^2$
has a point of self-intersection.
\end{proposition}

{\bf Proof.} Suppose by contradiction,
that $p(\gamma)$ has no self-intersection on $S^2$.
Denote by $\Omega_{\gamma}$ a disk on
$S^2$ bounded by $\gamma$.
Subtracting the first equation of \ref{twoineq}
from the second we get
\begin{equation}\label{congruence}
0\equiv\pi+\int_{\Omega_\gamma}(2-\frac{K}{2})\;{\rm dS}
\;\; mod \;\; l(S^1).
\end{equation}
One the other hand using $K\le 4$ with $\int_{S^2}\frac{K}{2} \;{\rm dS}=2\pi$,
 we get
\begin{equation}\label{long}
\pi\le \pi+\int_{\Omega_\gamma}(2-\frac{K}{2})\;{\rm dS}=
-\pi+2S(\Omega_{\gamma})+\int_{S^2\setminus \Omega_\gamma}\frac{K}{2}
{\rm dS}<2Area(S^2)-\pi,
\end{equation}
these inequalities contradict Equation (\ref{congruence})
since $l(S^1)=2Area(S^2)$ by \ref{2area}.

\hfill $\square$

\begin{remark} \label{slefinSK}

The statement of Proposition \ref{selfint}
holds as well if $S^3$ is isometric to the unit sphere of
a regular $PK$ cone, such that the conical
angles of the quotient $S^3_K/S^1$ are greater than $2\pi$.
The proof goes as above, one only needs to justify the equality
inside Equation (\ref{long}). To do this we follow the recipe of
Proposition \ref{twoineq} 2), namely for each conical point $x$
on $p(\gamma)$ we split the contribution of the curvature
$\delta$-function $\delta(x)$ as $\delta_{in}(x)+\delta_{out}(x)$,
with $\delta_{in}$ contributing to $\Omega_{\gamma}$, while
$\delta_{out}$ contributing to $S^2\setminus \Omega_{\gamma}$.

\end{remark}

Let $\phi$ be a closed curve on a real surface,
$\phi:S^{1}\to M^{2}$. A piece $\phi[t_{1},t_{2}]$ of the
curve $\phi$ is a {\it loop} if $\phi(t_1)=\phi(t_2)$.
The following lemma is trivial.

\begin{lemma} \label{leaves} A closed  curve with
self-intersections on a real surface has
at least two loops $\phi[t_1,t_2]$ and $\phi[t_3,t_4]$ that have
no overlap, i.e., $(t_1,t_2)\cap(t_3,t_4)=\emptyset$.
\end{lemma}

Let us now formulate and proof the smooth analog of
Theorem \ref{CAT1S3}.
\begin{theorem} \label{smoothCAT1S3}
Consider a regular Sasakian metric on $S^3$
such that  the quotient $S^3/S^1$ is $\rm CAT(4)$.
Then $S^3$ is $\rm CAT(1)$.
\end{theorem}
{\bf Proof.} By Proposition \ref{threestandard}  the sectional
curvature of $S^3$ is bounded from above by $1$,
so to prove that $S^3$ is $\rm CAT(1)$
we need to show that the length
of any closed geodesic $\gamma$  on $S^3$ is at least $2\pi$.

By Proposition \ref{selfint} the projection $p(\gamma)$ of $\gamma$ to
$S^3/S^1$ has at least one self-intersection. This means that $p(\gamma)$
contains at least two loops. By Theorem \ref{piloop},
each loop is no shorter than the complete
circle of the same geodesic curvature of the sphere of curvature $4$.
Hence by Proposition \ref{Bryant} the length of a piece of $\gamma$
that projects to a loop is at least $\pi$.
We deduce $l(\gamma)\ge 2\pi$.

\hfill $\square$

{\bf Proof of Theorem \ref{CAT1S3}.} This proof repeats
the proof of \ref {smoothCAT1S3}  with the following changes.
One needs to adjust Proposition \ref{selfint}
according to Remark \ref{slefinSK}. In order to apply (comparison)
Theorem \ref{piloop} to a loop of $p(\gamma)$
we need to show that the geodesic curvature
of $p(\gamma)$ at conical points of the surface (where
$p(\gamma)$ is not smooth) is not larger than at a smooth point
of the curve. This follows directly from definition
of geodesic curvature on page 70 of \cite{AB}.
Finally,
instead of Proposition \ref{Bryant} we use Lemma \ref{geonSK3}.

\hfill $\square$

\begin{remark} It would be interesting to show that these
results hold in higher dimensions,
for example to prove that a regular Sasakian
metric on $S^{2n+1}$ is $\rm CAT(1)$ if the metric on the
quotient $\mathbb CP^n\cong S^{2n+1}/S^1$ is $\rm CAT(4)$.

\end{remark}

\subsection{$\rm CAT(4)$ covers of positively curved $2$-spheres}

In this subsection we explain that a curvature $4$ sphere
with conical singularities of angles {\it less} than $2\pi$
admits ramified covers such that the metric pulled back to the cover
is $\rm CAT(4)$ (Corollary \ref{good1}).  First, we will need a result showing
that conical points on such a sphere are distributed "evenly".

\begin{proposition} \label{good} Let $S^2$ be a sphere
with a metric of curvature $4$ with
conical points $x_1,...,x_k$, $k>2$
of angles less than $2\pi$. Then
for any point $x\in S^2$ there exists  $x_i$
such that $d(x,x_i)<\frac{\pi}{4}$ ($x_i$ coincides with $x$
if $x$ is conical itself).
\end{proposition}

The proof of the proposition uses
a theorem of Toponogov and a lemma.

\begin{theorem}{\bf (Toponogov, \cite{To}.)}
Let $g$ be a metric on $S^2$ with Gaussian
curvature $K\ge 4$. Then any simple closed geodesic
on $S^2$ has length at most $\pi$. If one of such geodesics
has length $\pi$, then the metric $g$ is of curvature $4$.
\end{theorem}

\begin{lemma} \label{halfsphere}
For a sphere $S^2$ of curvature $4$
with conical points and a smooth point $x$ on it
the following statements hold:

1) Suppose that there is a geodesic loop $\gamma$
of length at most $\frac{\pi}{2}$ based
at $x$. Then there is a conical point
on $S^2$ of distance at most $\frac{\pi}{4}$ from $x$.

2) Suppose that the distance from
$x$ to all conical points is larger than $d$, and $d>\frac{\pi}{4}$.
Then $S^2$ contains an isometric copy of a spherical
disc of radius $d$ with the center at $x$.
\end{lemma}

{\bf Proof.}
1) The geodesic loop $\gamma$ cuts a neighborhood of $x$
into two sectors and the angle of one of the sectors
is at most $\pi$, denote this angle by $\alpha$.
Let $ABC$ be a geodesic triangle on the standard sphere of curvature $4$
such the length of $AB$ is $|\gamma|$ and the angles at $A$ and
$B$ are $\frac{\alpha}{2}$. Glue by isometry $CA$ with $CB$
identifying $B$ with $C$ so that we get a figure that we call a
{\it cup}. The {\it cup} has one conical point in its interior
(the point corresponding to
$C$), and by construction its
boundary is a geodesic loop of length $|\gamma|$ with angle
$\alpha$ at its base point. Notice that in  the triangle
$ABC$ the distance from  any point  to the set $\{A,B\}$
is at most $\frac{\pi}{4}$ ($ABC$ is contained
in a equilateral triangle with side $\frac{\pi}{4}$).
Hence  the distance from any point of
the {\it cup} to the angle on its boundary is at most
$\frac{\pi}{4}$.

Consider now a local isometry from a neighborhood of the boundary
of the {\it cup} to a neighborhood of $\gamma$. Suppose by contradiction
that all conical points on $S^2$ are on distance
at least $\frac{\pi}{4}$ to $x$. Then this local isometry can
be extended to the local isometry from the {\it cup}
minus its conical point. Hence it extends by continuity
to the conical point, whose image is a conical point on $S^2$
on the distance at most $\frac{\pi}{4}$ from $x$.

2) Since there are no conical points on $S^2$ on distance
less than $\frac{\pi}{4}$ from $x$, it follows
from 1) that the radius of injectivity at $x$
is at least $\frac{\pi}{4}$. Hence $S^2$ contains an isometric
copy of the discs of radius $\frac{\pi}{4}$ centered at $x$.
Consider now the family of disks centered at $x$ of radius growing
from $\frac{\pi}{4}$ to $d$. The boundaries of these disks are smooth
(they don't contain conical points) and  concave, hence
they can not acquire a self-intersection as the radius grows.
Therefore  they are all embedded in the sphere.

\hfill $\square$

{\bf Proof of Proposition \ref{good}.}
Suppose first by contradiction that there
is a point on $S^2$ on distance
more than $\frac{\pi}{4}$ from each conical point.
Then by Lemma \ref{halfsphere} 2) the sphere contains an isometric copy
of a radius $\frac{\pi}{4}$ disk.
Its boundary is a simple closed geodesic of
length $\pi$. This contradicts Toponogov's theorem,
since we can smoothen the metric
near conical points so that the curvature is still at least
$4$ and the geodesic is kept intact.

To finish the proof we need to rule out the case when
the distance from a point $x$ on $S^2$ to all conical points
is a least $\frac{\pi}{4}$ and for some $i$
$d(x,x_i)=\frac{\pi}{4}$. By arguments similar to those used
in the proof of Lemma \ref{halfsphere} one can check in this case
that  $S^2$ has exactly two conical points and
$x$ lies on the closed geodesic of length less than $\pi$
that splits $S^2$ into two equal parts.
But by the assumptions the number of conical points on $S^2$
is at least $3$.

\hfill $\square$

\begin{corollary} \label{good1} Let $S^2$ be a sphere
with a metric of curvature $4$ with
conical points $x_1,...,x_k$, $k>2$
of angles less than $2\pi$. Then we have:

1) There exists a geodesic triangulation of $S^2$ by convex triangles with vertices at  $x_1,...,x_k$.

2) For sufficiently  large multiplicities $n_1,...,n_k$
the universal orbifold cover of  the orbifold
$(S^2, x_1,n_1,...,x_k, n_k)$ is a $\rm CAT(4)$ space.

3) For any cover from 2) there exists $C>0$ such that
simple closed curves $\gamma$
on the cover satisfy the isoperimetric
inequality $Area(\gamma)< C\cdot l(\gamma)^2$.

\end{corollary}

{\bf Proof of Corollary \ref{good1}.1)}
This proof
is almost identical to the proof of Proposition 3.1
in \cite{Th} claiming the existence of a triangulation
for a sphere with a flat metric and with conical points
of angles less than $2\pi$. We will just recall the main idea.

One shows first that there exists a canonical decomposition of
$S^2$ into a collection of convex polygons whose edges are geodesics
joining conical points of $S^2$. To construct the polygons
one considers all (finite number)
maps from round disks to $S^2$ that are locally
isometric on the interior and such that at least $3$ conical
points are contained in the image of the boundary of the disk.
The polygon is obtained as the convex hull of all conical points
on the disk boundary. Once polygons are constructed each of them
can be decomposed into triangles. Notice that we use
here Proposition \ref{good} to assure the all considered
disks are of radius less than $\frac{\pi}{4}$.

\hfill $\square$

The proof of \ref{good1} 2) is slightly less elementary
than one would want it to be but at least it avoids
calculations.
We will use a couple of results of Bowditch and a simple lemma.

\begin{theorem} \label{Bowditch}
{\bf(Bowditch \cite{Bo}.)} In a
complete locally compact locally $\rm CAT(1)$ space a closed geodesic
cannot be freely homotoped to a point through
rectifiable curves of length strictly less than $2\pi$.
\end{theorem}

\begin{lemma} \label{shortening}
Any piecewise smooth closed
curve on a compact constant curvature
surface with conical points can be freely homotoped to
a closed geodesic or a point via rectifiable curves
of decreasing length.
\end{lemma}
This last lemma is a very special case of
Theorem 3.1.6 from \cite{Bo}.

\begin{lemma} \label{liptriangle}
For a convex
spherical triangle $\Delta_4$  of curvature $4$
and its comparison triangle $\Delta_{\kappa}$
of curvature $\kappa\le 4$
there exists a bi-Lipschitz map $\Delta_4\to \Delta_\kappa$
that is an isometry on the boundary. Moreover the bi-Lipschitz
constant is bounded by a continuous function $c(\Delta_4,\kappa)$
with $c(\Delta_4,4)=1$.

\end{lemma}
{\bf Proof.} Let $O_4$ and $O_{\kappa}$ be the centers
of circles inscribed in $\Delta_4$ and $\Delta_{\kappa}$
correspondingly. Consider the map sending isometrically the boundary
of $\Delta_4$ to the boundary of $\Delta_{\kappa}$
and sending geodesic segments
through $O_4$ to geodesic segments through $O_{\kappa}$
so that the metric restricted to each
geodesic segment is multiplied by a constant. A calculation
shows that this map is bi-Lipschtiz and the bi-Lipschtiz
constant tends to $1$ as $\kappa$ tends to $4$.

\hfill $\square$

{\bf Proof of  Corollary \ref{good1}.2)}
Consider the triangulation of $S^2$ by convex triangles given by
\ref{good} 1).
Replacing each convex triangle by a comparison triangle of
curvature $0\le \kappa\le 4$  we get a family of spheres
$S^2_{\kappa}$ of curvature $\kappa$ with conical points,
$S^2_4$ being the original sphere.
The cone angles of $S^2_{\kappa}$ decrease with $\kappa$,
and we denote by $\alpha_{i0}$ the cone angle of $S^2_0$ at $x_i$.
The sphere $S^2_0$ is of curvature $0$ and for $n_1,...,n_k$
such that $n_i\alpha_{i0}\ge 2\pi$
the orbifold universal cover $\widetilde S^2_{0}$
of $(S^2_{0},n_i,x_i)$ is $\rm CAT(0)$.
It follows that $\widetilde S^2_{0}$  does not contain closed
geodesics. We will deduce from this that
$\widetilde S^2_{4}$  does not
contain closed geodesics of length $\pi$
or less, and so it is $\rm CAT(4)$.

Suppose by contradiction that on $\widetilde S^2_4$ there is a
geodesic of $\gamma$
length $\pi-\varepsilon$. Denote by $T$ the subset of
$[0,4]$ of  $\kappa$ such that $\widetilde S^2_{\kappa}$
contains a geodesic of length at most $\pi-\varepsilon$.
By standard arguments $T$ is closed
and we will show that it is also open. This will
imply that $T=[0,4]$ and will lead to a contradiction
since $0\notin T$. To show that $T$ is open
we need to prove that for each $\kappa\in T$
a neighborhood of $\kappa$ belongs to $T$, we will
do it for $\kappa=4$, for other points the reasoning is
identical.

Chose $\delta$ so that for each triangle $\Delta_4$
from the convex triangulation of $\widetilde S^2_4$ we have
$c(\Delta_4,\kappa)<\sqrt{\frac{\pi}{\pi-\varepsilon}}$ for
$\kappa\in [4-\delta,4]$.
Let us show that  for $\kappa\in [4-\delta,4]$ $\widetilde S^2_{\kappa}$ contains
a geodesic of length less than $\pi-\varepsilon$.
 By Lemma \ref{liptriangle} for
$\kappa\ge 4-\delta$ there exists a bi-Lipschitz map $F_{\kappa}$
from $\widetilde S^2_4$ to $\widetilde S^2_{\kappa}$ with bi-Lipschitz
constant less than $\sqrt{\frac{\pi}{\pi-\varepsilon}}$.
Notice that the curve $F_{\kappa}(\gamma)$
can not be contracted to a point on $\widetilde S^2_{\kappa}$
by a length decreasing homotopy.
Indeed, during such a homotopy the length of the curve would
be at most $l(F_{\kappa}(\gamma))\le \sqrt{\pi(\pi-\varepsilon)}$.
Hence the pull back of such a homotopy to $S^2_4$ would
contract  $\gamma$
to a point keeping its length less than $\pi$, which
contradicts Theorem \ref{Bowditch}.

Denote by $\gamma'$ the piecewise geodesic curve on
$\widetilde S^2_{\kappa}$ obtained from $F_{\kappa}(\gamma)$
by replacing each interval of $F_{\kappa}(\gamma)$
with the ends at the boundary of the triangulation by
a geodesic segment. The curve  $F_{\kappa}(\gamma)$
can be deformed to $\gamma'$ by a length decreasing
deformation. By Lemma \ref{shortening} $\gamma'$
can be homotoped further to a shorter geodesic $\gamma''$.
Finally, it follows from the  $\rm CAT(\kappa)$ inequality
(Definition \ref{catdefinition}) that
$l(\gamma')\le l(\gamma)$. Hence as promised we constructed
on $\widetilde S^2_{\kappa}$ a closed geodesic $\gamma''$ with
$l(\gamma'')<\pi-\varepsilon$. This finishes the proof that
$\widetilde S^2_4$ is $\rm CAT(4)$.

\hfill $\square$

{\bf Proof of  Corollary \ref{good1}.3)}
In order to prove that  curves on
$\widetilde S^2_4$ satisfy the isoperimetric
inequality recall that there is a bi-Lipschitz map from
$\widetilde S^2_4$ to $\widetilde S^2_0$, and
the last space is $\rm CAT(0)$, hence curves
on it satisfy the desired isoperimetric inequality.

\hfill $\square$

\subsection{Proof of Theorem \ref{COVERCAT1} and its refinements}

Let $C_{K}^{4}$ be a regular spherical $PK$ cone. According
to Theorem 1.7 \cite{P} the cone is canonically (up to a scale)
biholomorphic to $\mathbb C^2$
and its singular locus is defined by a collection of linear
equations $L_{1}=0,...,L_{k}=0$ in $\mathbb C^2$.
Moreover, the multiplication of $\mathbb C^2$
by unitary complex numbers  generates
a free isometric $S^1$-action on $C_K^4$.
For $n>1$ consider the ramified (Galois)
covering of $C_{K}^{4}$ of degree $n^{k}$ corresponding to the
extension of the field of rational functions on $\mathbb{C}^{2}$
by functions $L_{1}^{\frac{1}{n}},...,L_{k}^{\frac{1}{n}}$. The
corresponding $PK$ cone that covers $C_{K}^{4}$ is denoted by
$(C_{K}^{4})^{\frac{1}{n}}$, the unit sphere of this cone
is denoted by $(S_K^3)^{\frac{1}{n}}$.  One can check that the
isometric $S^1$-action on $C^4_K$ generated a {\it free}
isometric action on
 $(C_K^4)^{\frac{1}{n}}$
(and as well on its unit sphere $(S_K^3)^{\frac{1}{n}}$).

\begin{proposition}\label{negativecone}
Let $C_{K}^{4}$ be a {\it positively curved} regular spherical $PK$ cone.
For sufficiently large $n$
every closed geodesic $\gamma$ in the unit sphere
$(S_K^3)^{\frac{1}{n}}$ of the cone $(C_{K}^{4})^{\frac{1}{n}}$
with contractible projection $p(\gamma)$
on  $(S_K^3)^{\frac{1}{n}}/S^1$ has length at least $2\pi$.
\end{proposition}

{\bf Proof.} Denote by $M_n^2$ the universal
cover of $(S_K^3)^{\frac{1}{n}}/S^1$ and
by $M_n^3$ the corresponding cover of $(S_K^3)^{\frac{1}{n}}$.
Then any geodesic $\gamma$ in $(S_K^3)^{\frac{1}{n}}$ with contractible
projection on $(S_K^3)^{\frac{1}{n}}/S^1$ admits a lift to $M_n^3$.
So it is sufficient
to show that for $n$ large enough any closed geodesic $\gamma'$ in
$M_n^3$ is longer than $2\pi$.

Note that the  cover $M_n^2$ is isomorphic to the orbifold
universal cover of the sphere $S_K^3/S^1$ with orbifold
structure of multiplicity $n$ at all its conical points.
Hence by Corollary  \ref{good1}.2) $M_n^2$  is $\rm CAT(4)$  for
$n$ large enough. If $p(\gamma')$ has a self-intersection
then it has at least two loops  each contributing
at least $\pi$ to the length of $\gamma$ and so
$l(\gamma')\ge 2\pi$ (see Theorem \ref{piloop} and
the proof of Theorem \ref{smoothCAT1S3}).
Recall now
Equation (\ref{modl(S)}) and (\ref{gaussbonnet})
$$\int_{\gamma'}w_0\equiv
-2Area (\Omega_{\gamma'})\;\; {\rm mod}\;\; l(S^1),
\;\;\;
\int_{\gamma'}w_0=
\frac{1}{2}(2\pi-\int_{\Omega_{\gamma'}}K\;{\rm dS}).$$

Using $K\le 4$ we get
\begin{equation}\label{convexmust}
\int_{\gamma'}w_0=
\frac{1}{2}(2\pi-\int_{\Omega_{\gamma'}}K\;{\rm dS})\ge
\pi-\int_{\Omega_{\gamma'}}2{\rm dS}>-2Area (\Omega_{\gamma'}).
\end{equation}

Suppose now that $l(\gamma')<2\pi$.
Then the isoperimetric inequality
(Corollary \ref{good1} 3)) tells us
$Area(\Omega_{\gamma'})<C(4\pi^2)$.
At the same time we can chose $n$ such that
$l(S^1)>2\pi+2C(4\pi^2)$.
For such $n$ Equation (\ref{modl(S)}) with
Inequality (\ref{convexmust}) imply
$\int_{\gamma'}w_0>2\pi$. Hence $l(\gamma')>2\pi$
and we get a contradiction.

\hfill $\square$

Proposition \ref{negativecone} should be compared with the following lemma.

\begin{lemma}\label{noncat}
 Suppose that $C^4_K$
is not a direct product of two
$2$-dimensional cones, and let $\alpha_{min}$
be the minimal conical angle at
the singular locus of $C^4_K$. If
$\alpha_{min}\cdot[\frac{n}{2}]\ge \pi$ then
$(C^4_K)^\frac{1}{n}$
is not locally $\rm CAT(0)$.
\end{lemma}

{\bf Proof.}
Indeed, if $C^4_K$ is not a direct product
then any two conical points $x$ and $y$ on the sphere
$S^3_K/S^1$ can be joined by a geodesic
segment $[xy]$ with $l([xy])<\frac{\pi}{2}$.
Consider in $S^3_K$ a geodesic segment  $[x'y']$
that joins singular Hopf fibers over $x$ and $y$
and projects to $[xy]$. This segment is orthogonal
to the $S^1$ fibers and has the same length as $[xy]$.
Consider now on the complement to singular circles
in $S^3_K$ the following path contained in an
$\varepsilon$-neighborhood
of $[x'y']$. The path starts in the middle of $[x'y']$,
goes to the singular circle  containing $x'$,
enlaces it $[\frac{n}{2}]$ times, goes to $y'$,
enlaces the singular circle $[\frac{n}{2}]$ times and
comes back to the middle of $[x'y']$, after this the path repeats
itself again, but enlaces the circles in the opposite
direction. If follows from the construction of
$(C^4_K)^{\frac{1}{n}}$ that this path lifts
to a closed path on $(S^3_K)^{\frac{1}{n}}$.
Moreover in the case when
the conical angles at $x$ and $y$ satisfy
$\alpha(x)[\frac{n}{2}]>\pi$ and
$\alpha(y)[\frac{n}{2}]>\pi$ there is a closed
geodesic in $(S^3_K)^{\frac{1}{n}}$ composed
of $4$ preimages of $[x'y']$, contained in
the $\varepsilon$-neighborhood of the lifted path.

\hfill $\square$

\begin{remark} Notice that the cone $(C^4_K)^\frac{1}{n}$
is not $\rm CAT(0)$ if $\alpha_{min}\cdot n<2\pi$, so it is reasonable
to guess that it can be $\rm CAT(0)$ only
if $C^4_K$ is a direct product of two  $2$-cones. Lemma
\ref{noncat} also
explains why the degree $n^2$ cover $\mathbb C^2\to \mathbb C^2$ with
order $n$ branching at two non-orthogonal lines is not $\rm CAT(0)$.
\end{remark}

The following Proposition provides an explicit construction
for Theorem \ref{COVERCAT1}.

\begin{proposition}
\label{orbipoint} For every cone $(C_K^4)^{\frac{1}{n}}$
satisfying the conditions of Proposition \ref{negativecone} there
exists a cone $\widetilde C_K^4$ satisfying the following properties:

1) There is an action of a finite
group $G$ on  $\widetilde{C_K^4}$ such that
$(C_K^4(x_i))^{\frac{1}{n}}=\widetilde{C_K^4}/G$; the action of
$G$ is free outside of the origin of $\widetilde{C_K^4}$.

2) Closed geodesics in the unit sphere $\widetilde{S_K^3}$
are longer than $2\pi$.

3) The map of $\widetilde{C_K^4}$ blown up at the
origin to $(C_K^4)^{\frac{1}{n}}$ blown up at the origin is
a non-ramified covering.
\end{proposition}

{\bf Proof.} Consider the Riemann surface
$M_n=(S_K^3)^{\frac{1}{n}}/S^1$. Since $\pi_1(M_n)$
is a residually finite group (see \cite{He} for an elementary proof)
there exists a Galois covering $\widehat{M_n}\to M_n$ with finite Galois group $G$
such that all non-contractible geodesics on $\widehat{M_n}$
are longer than $2\pi$.

Denote by $(C_K^4)_{\circ}^{\frac{1}{n}}$ the cone $(C_K^4)^{\frac{1}{n}}$
blown up at the origin. This is
the total space of a certain  line bundle over $M_n$ and so
$\pi_1((C_K^4)_{\circ}^{\frac{1}{n}})\cong \pi_1(M_n)$.
Take the Galois covering of
$(C_K^4)_{\circ}^{\frac{1}{n}}$ with the Galois group $G$.
Then the cone $\widetilde {C_K^4}$
obtained by contraction of the exceptional curve on the constructed
covering, satisfies the properties 1), 2), 3). Indeed, 2) is verified
because any geodesic $\gamma$ on $\widetilde {S_K^3}$ with contractible
$p(\gamma)$ is longer than $2\pi$ by Proposition \ref{negativecone};
if $p(\gamma)$ is non-contractible, then by construction $l(\gamma)\ge
l(p(\gamma))>2\pi$. The properties 1) and 3) follows automatically
from the construction of  $\widetilde {C_K^4}$.

\hfill $\square$

{\bf Proof of Theorem \ref{COVERCAT1}.} Chose $n$ such that the cone
$(C_K^4)^{\frac{1}{n}}$ satisfies the conditions of Proposition
\ref{negativecone} and take its cover provided by Proposition
\ref{orbipoint}. The unit sphere of the obtained cone
is the desired $3$-manifold.

\hfill $\square$

\begin{remark} In Theorem \ref{COVERCAT1} one can let conical angles tend to $2\pi$. As a result
we get the statement that for each collection $L_1,...,L_k$ of complex lines in a flat
$\mathbb C^2$ going trough $0$ one can construct a locally $\rm CAT(0)$ ramified cover of $\mathbb C^2$ with branchings at $L_i$, provided the following condition holds.
For the points $p_1,...,p_k$ on the projectivization $\mathbb CP^1$ of $\mathbb C^2$
corresponding to $L_1,...,L_k$ there is no geodesic $S^1$ in
$\mathbb CP^1$ such that $p_1,...,p_k$ are contained
in one connected component of  $\mathbb CP^1\setminus S^1$.
\end{remark}

\section{Construction of complex surfaces}

\subsection{Proof of Theorem \ref{ramCP2}}

In order to prove Theorem \ref{ramCP2}
we need to use the notion of {\it rigid orbispaces}
due to Haefliger, defined in \cite{GH}, chapter 11.
This is a natural  generalization of the notion of orbifolds
to the category of locally compact topological spaces.
Rigid   orbispace is defined by replacing the model
for the orbifold charts by a locally compact space
with a rigid action of a finite group,
i.e. one for which points with trivial isotropy are dense.
The precise definition can be found in \cite{GH} and \cite{CD}.
We say that a geodesic metric space $X$ with a
structure of a rigid orbispace is
{\it non-positively curved} if the
metric induced on each model chart is
locally $\rm CAT(0)$.

The following theorem is proven in Chapter 11 of \cite{GH}.
\begin{theorem}\label{thm:orbi}
Any non-positively curved  rigid orbispace $X$
admits a universal cover, i.e., there exists a simply connected
topological space with an action of a group $G$ such that
$Y/G=X$, $stab(y)$ is finite for any $y\in Y$, and the
orbispace structure induced on $X$ by the quotient coincides with
the original one.
\end{theorem}

{\bf Proof of Theorem \ref{ramCP2}.}
Let $x_1,...,x_m$ be the multiple
points of the arrangement. Let $C_K^4(x_1),...,C_K^4(x_m)$ be
the tangent cones to $\mathbb {C}P^2$ at the points $x_1,...,x_m$.
Chose such $N'$ that for any $n\ge N'$ the cones
$(C_K^4(x_1))^{\frac{1}{n}},...,(C_K^4(x_m))^{\frac{1}{n}}$
satisfy the condition of Proposition  \ref{negativecone}.
Put
$$N=\max\Bigl(N', \max_j \frac{1}{\beta_j}\Bigr).$$

Pullback the
metric from $\mathbb CP^2$ to the complex surface $S(n,L_{j})$.
If $n\ge N$ the  conical angles at
singular curves on $S(n,L_{j})$ are greater than $2\pi$,
and so the metric has
non-positive curvature on the complement to the
pre-images of the points $x_1,...,x_m$.
At the same time each point $\widetilde x_i$ in the preimage of  $x_i$
has a neighborhood isometric to a neighborhood of the origin in
the cone  $(C_K^4(x_i))^{\frac{1}{n}}$.
Using Proposition \ref{orbipoint} we can
represent each cone $(C_K^4(x_i))^{\frac{1}{n}}$
as a quotient of a locally $\rm CAT(0)$ cone by a finite group
(where only the origin has a nontrivial stabiliser).
This defines  on the complex
surface $S(n,L_j)$ a structure of a
rigid non-positively curved orbispace.

By Theorem \ref{thm:orbi} the complex surface
$S(n,L_j)$ with  the constructed rigid orbispace structure admits a
universal orbispace cover $\widetilde S$. The pullback metric on
$\widetilde S$
has non-positive curvature, so by Theorem \ref{kpi1}
the complex surface $\widetilde S_{\circ}$ obtained from
$\widetilde S$ by the
blowup of its  complex singularities is of the type $K(\pi,1)$.
At the same time by condition 3) of Proposition
\ref{orbipoint} the complex surface $\widetilde S_{\circ}$
is a non-ramified cover of the blow up of $S(n,L_j)$
at its complex singularities.
Hence  $S(n,L_j)$ is of the type $K(\pi,1)$ as well.

 \hfill $\square$

\subsection{\label{orbiresults} Orbifold structures on complex surfaces}

Consider a smooth complex surface $S$ with a divisor
$D=\sum_{i=1}^n D_i$ such that the divisors $D_i$ are
smooth, irreducible and pairwise transversal. Suppose we want
to introduce on $S$ the structure of an orbifold
so that for points of $D_i$ that are not multiple in $D$
the local fundamental group is
$\mathbb Z_{b_i}$, $b_i>1$. This imposes conditions
on the multiple points of $D$ and on the numbers $b_i$
summarised in the next proposition
(see for example \cite{Y} Section 11, or \cite{U}).

\begin{proposition} \label{orbiglobal}
The complex surface $(S,D)$ admits an orbifold
structure with local fundamental group
$\mathbb Z_{b_i}$ at generic points of $D_i$ if all multiple points of $D$ are at most triple, and at each triple point $\frac{1}{b_j}+\frac{1}{b_k}
+\frac{1}{b_l}>1$ for the incoming multiplicities.
The order of the local
fundamental group at the triple point is
$4(\frac{1}{b_j}+\frac{1}{b_k}+\frac{1}{b_l}-1)^{-2}$.
\end{proposition}

This proposition relies on the list of Shephard and Todd
of complex reflection groups acting on $\mathbb C^2$.
For each triple point of the divisor infinitesimally one has to get
an action of $G$ on $\mathbb C^2$ such that
$\mathbb C^2/G\cong \mathbb C^2$ and moreover the quotient
map  $\mathbb C^2\to \mathbb C^2/G$
is a branched cover that ramifies along three lines
in $\mathbb C^2/G\cong \mathbb C^2$. For each triple
of branching multiplicities $p,q,r$ satisfying
$\frac{1}{p}+\frac{1}{q}+\frac{1}{r}>1$
there is a unique group in the list of Shephard and Todd.
The group is called $G(2s,2,2)$ in the case $(p=s,q=r=2)$, and has numbers
$(7)$, $(11)$, and $(19)$ for triples $(2,3,3)$, $(2,3,4)$, and $(2,3,5)$
correspondingly. These groups act on the projectivisation
$\mathbb CP^1$ of $\mathbb C^2$ and the action factors
through the action of a dihedral group for $G(2s,2,2)$,
and tetrahedral, octahedral, and dodecahedral groups
in the other $3$ cases.

{\bf Orbifold Chern numbers and Miyoka-Yau.}
Recall that according to a result of Miyaoka and Yau
for a complex surface $S$ of general type the inequality on
Chern numbers holds $c_1^2(S)\le 3c_2(S)$ and
the equality is attained if and only if the universal
cover of $S$ is the unit ball $|z_1|^2+|z_2|^2<1$.
The analogous statement  for orbifolds is contained
in \cite{KNS}. We will only state a very special
case of this theorem when $S\cong \mathbb CP^2$ and
the divisor $D=D_1+...+D_n$ satisfies the conditions of Proposition
\ref{orbiglobal}. Denote by $d_i$ the degree of $D_i$
and for each singular point $p$ of $D$ denote
by $\beta_D(p)$ the order of the local group at $p$.
Then the Chern numbers of the orbifold $(\mathbb CP^2,D,b)$
are defined as follows.

\begin{equation} \label{orbic1}
c_1^2(\mathbb CP^2,D,b)=(-3+\sum_{1\le i\le n}d_i(1-b_i^{-1}))^2
\end{equation}
\begin{equation}\label{orbic2}
e(\mathbb CP^2,D,b)=3-
\sum_{1\le i\le n}(1-b_i^{-1})e(D_i\setminus sing(D))-
\sum_{sing(D)}(1-\beta_D(p)^{-1})
\end{equation}
The following theorem is a special
case of a theorem from \cite{KNS}.
\begin{theorem} \label{KNS}
Let $(\mathbb CP^2,D,b)$
be an orbifold that is not a global quotient of $\mathbb CP^2$.
Then $c_1^2(\mathbb CP^2, D,b)\le 3e(\mathbb CP^2,D,b)$,
the equality holding if and only if $(\mathbb CP^2,D,b)$
is uniformized by the complex ball.
\end{theorem}

\begin{corollary}\label{plball} The orbifold $(\mathbb CP^2,A^0_3(3),2)$
is uniformized by the complex ball.
\end{corollary}

{\bf Proof.} The arrangement $A_3^0(3)$ consist of
$9$ lines each of which intersect others in $4$ points.
Each multiple point $p$ is triple with $\beta_D(p)=16$,
and there are $12$ of them.
Applying equations (\ref{orbic1}) and
(\ref{orbic2}) we get
$$c_1^2(\mathbb CP^2,A_3^0(3),2)=\frac{9}{4}=3e(\mathbb CP^2,A_3^0(3),2),$$
and since $(\mathbb CP^2,A_3^0(3),2)$ is not uniformized
by $\mathbb CP^2$ the proof is finished by
Theorem \ref{KNS}.
\hfill $\square$

\subsection{\label{threearrangements} Proof of Theorem \ref{CPorbi}}

Let us first give a description of the arrangements $A_1(6)$, $A_1(7)$, and $A^0_3(3)$.

1) The {\it complete quadrilateral} $A_1(6)$ consists of $6$
lines on $\mathbb CP^2$  that pass through $4$ generic points.
The corresponding $PK$ metric on $\mathbb CP^2$ has conical
angles $\pi$ at all lines.

2) $A_1(7)$ consists of lines $x=\pm z$, $y=\pm z$, $x=\pm y$,
and $z=0$.  The  metric has angles
$\pi$ at  lines $x=\pm z$, $y=\pm z$, and angles $\frac{4\pi}{3}$
at the rest.

3) $A^0_3(3)$ is the arrangement of $9$ lines
given by the equations $x^3-y^3=0$, $y^3-z^3=0$ and $z^3-y^3=0$
with angles $\frac{4\pi}{3}$ at all of them.

Theorem \ref{CPorbi} will be deduced from Theorem \ref{CAT1S3}
and a special case of a lemma from \cite{Gr1}.
Call a convex triangle on $S^2$ of curvature
$4$ {\it large} if the distance from each
vertex to the opposite side is a least
$\frac{\pi}{4}$.

\begin{lemma} \label{grompi4}
( {\bf 4.2.E,  Remark b, \cite{Gr1}.})
Consider a real surface of curvature
$4$ with conical singularities  admitting a
triangulation by large triangles of valence
at least four.
If any two conical points are joined by
at most one edge and any cycle consisting of
$3$ edges borders a face of the triangulation then
the surface is $\rm CAT(4)$.
\end{lemma}

{\bf Sketch of a proof.}
Let $\gamma$ be a closed geodesic on the surface. If $\gamma$ is
composed entirely of the edges of the triangulation then it should
contain at least four edges, hence $l(\gamma)\ge \pi$. Otherwise
$\gamma$ should intersect an edge $e$ of the triangulation. Let
$\Delta_1$ and $\Delta_2$ be the triangles adjacent to $e$ and
let $A_1$, $A_2$ be the vertices of $\Delta_1$ and $\Delta_2$
opposite to $e$. Let $L(A_1)$ and $L(A_2)$ be the links of
$A_1$ and $A_2$ consisting of all triangles containing $A_1$
and $A_2$ correspondingly.
From the condition of the lemma it follows that the interiors
of $L(A_1)$ and $L(A_2)$ don't intersect. The main observation
is that if a closed geodesic intersects a link of a vertex
on the surface then the length of the intersection is
at least $\frac{\pi}{2}$.

\hfill $\square$

{\bf Proof of Theorem \ref{CPorbi}.}
For the arrangements $A_1(6)$, $A_1(7)$ and $A^0_3(3)$ in  \cite{P}
were constructed $PK$  metrics on $\mathbb CP^2$
such that  all conical angles
along the lines of the arrangements are at least $\pi$.
In the case of $A_1(6)$ the conical angles at
all six lines are $\pi$,
for $A_1(7)$ four angles are $\pi$ and three angles are
$\frac{4\pi}{3}$, for $A^0_3(3)$ the angles at all
nine lines are $\frac{4\pi}{3}$.
Since the orbifold structure along each line
is non trivial the orbifold metrics have non-positive
curvature in the complement of the multiple points of the arrangements.

Let $x$ be a triple  point of $A_1(6)$, $A_1(7)$, or $A^0_3(3)$.
In order to prove that the orbifold has non-positive curvature at
$x$ we need to show  that the
orbifold cover of the unit sphere $S_K^3(x)$ is
$\rm CAT (1)$. By construction this cover is the unit
sphere of a regular spherical $PK$ cone, and so
by Theorem   \ref{CAT1S3}
it is sufficient to prove that the orbifold cover of  $S_K^3(x)/S^1$
is $\rm CAT (4)$.

The sphere $S_K^3(x)/S^1$ is a $2$-sphere of curvature
$4$ with three  conical points whose angles
are as follows: $(\pi,\pi, \pi)$ for $A_1(6)$,
$(\frac{4\pi}{3},\frac{4\pi}{3}, \pi)$ for $A_1(7)$, and
$(\frac{4\pi}{3},\frac{4\pi}{3},\frac{4\pi}{3})$ for $A^0_3(3)$.
In each case the sphere is obtained by gluing two copies of a
large triangle. Hence  the corresponding orbifold covers of $S_K^3(x)/S^1$
are triangulated by large triangles and it is straightforward
to check that these triangulations
satisfy the conditions of Lemma \ref{grompi4}.

\hfill $\square$

\begin{remark} \label{cat(-1)}The statement of Theorem  \ref{CPorbi} can be strengthen
in the case of arrangement $A^0_3(3)$. It follows from \cite{CHL}
that the Fubini-Studi metric on $\mathbb CP^2$ admits a one parameter deformation 
in the class of metrics of constant curvature with singularities at the arrangement $A^0_3(3)$;
the family is parameterized by the conical angle at the lines of  $A^0_3(3)$.
Both the flat metric that we used and the ball quotient metric (from Corollary \ref{plball}) belong to this family, corresponding conical angles being $\frac{4\pi}{3}$ and $\pi$.
Using the methods of the present article
one can show that each metric from the family with conical angles in 
the interval $(\pi,\frac{4\pi}{3})$ is
locally $\rm CAT(K)$, $K<0$ with respect to any orbifold structure from Theorem  \ref{CPorbi}.

\end{remark}

{\bf Locally $\rm CAT(0)$ polyhedral metrics on hyperbolic manifolds.}
Theorem \ref{CPorbi} can be used to construct a locally
$\rm CAT(0)$ polyhedral metric on a compact complex surface that is a smooth  quotient of the
unit two-dimensional complex ball. Namely,
we should apply Theorem \ref{CPorbi}
to the arrangement $A^0_3(3)$ with choices
$d_1=...=d_9=2$. By Corollary \ref{plball}
the orbifold universal cover of $\mathbb CP^2$ is the unit complex
ball $B^2$. So the orbifold fundamental group
is a co-compact lattice in $U(2,1)$, and taking a quotient of $B^2$
by its torsion free finite index subgroup we get an example.

By \cite{CDM} it is known that every real hyperbolic manifold
(of sectional curvature $-1$) admits a locally $\rm CAT(0)$
polyhedral metric, but the existence of such a metric
on a complex hyperbolic manifold is rather unexpected
(see the discussion on page 45 in \cite{DM}).
In \cite{Gr2} Gromov introduces
a class of {\it polyhedral
groups with $K<0$}, that are fundamental groups of
locally $\rm CAT(0)$ polyhedrons, satisfying
additional property of {\it non-flatness at each face} (page 176).
He makes a conjecture  (Section 7.A.IV, page 180)
that for every group of this
class any homomorphism  to a co-compact lattice in $U(n,1)$
should factor through a map to $\pi_1$ of a Riemann surface.
If the polyhedral
metrics that we construct on ball quotients admit a polyhedral decomposition,
such that each face is {\it not flat} in the sense of \cite{Gr2},
one would get a counter-example to the above conjecture.
It should not be difficult to check if we actually get a counterexample
but we have not yet done that.

\end{document}